\documentclass[12pt,leqno,twoside]{amsart}
\usepackage{amssymb}
\usepackage{latexsym}
\usepackage{verbatim}
\usepackage{epsfig}

\newtheorem{theorem}{Theorem}[section]
\newtheorem{corollary}[theorem]{Corollary}
\newtheorem{lemma}[theorem]{Lemma}
\newtheorem{proposition}[theorem]{Proposition}
\theoremstyle{definition}
\newtheorem{definition}[theorem]{Definition}
\theoremstyle{remark}
\newtheorem{remark}[theorem]{\sc Remark}
\theoremstyle{remark}

\theoremstyle{remark}
\newtheorem{example}[theorem]{\sc Example}
\theoremstyle{remark}
\newtheorem{note}[theorem]{\sc Note}
\theoremstyle{remark}

\theoremstyle{remark}

\renewcommand{\Box}{\square}    
\renewcommand{\Bbb}{\mathbb}
\newcommand{\cal}{\mathcal}

\newcommand{\Cone}{{\rm{Cone}}}
\renewcommand{\rhd}{{\rm{rhd\hspace{1.5pt}}}}
\newcommand{\h}{{\rm{ht}}}

\newcommand{\diff}{{\rm{diff}}}
\renewcommand{\int}{{\rm{int}}}

\newcommand{\Sing}{{\rm{Sing\hspace{1pt}}}}

\newcommand{\grad}{\mathop{\rm{grad}}\nolimits}

\newcommand{\fin}{\hspace*{\fill}$\Box$}
\newcommand{\m}{\setminus}
\renewcommand{\a}{\alpha}

\newcommand{\cC}{{\cal C}}
\newcommand{\cD}{{\cal D}}

\newcommand{\cS}{{\cal S}}

\newcommand{\cW}{{\cal W}}

\newcommand{\bC}{{\Bbb C}}

\newcommand{\bP}{{\Bbb P}}

\newcommand{\bX}{{\Bbb X}}
\newcommand{\bY}{{\Bbb Y}}

\newcommand{\bH}{{\Bbb H}}

\def\growarrow#1{
  \setbox1=\hbox{ $\scriptstyle #1$\ }

\mathop{\smash{\hbox to \wd1{\rightarrowfill}}
          \vphantom\rightarrow}\limits^{#1}}

\begin{document}

\title[Connectivity via nongeneric pencils]{Connectivity via nongeneric pencils}
\author{Mihai Tib\u ar}

\address{Math\' ematiques, UMR 8524 CNRS,
Universit\'e des Sciences et Tech. de Lille, \  59655 Villeneuve d'Ascq, France.}

\email{tibar@agat.univ-lille1.fr}

\thanks{This work was partially supported by the Newton Institute at Cambridge (Singularity Semester 2000).
}



\maketitle


\begin{abstract}
We use slicing by nongeneric pencils of hypersurfaces and prove a new theorem of Lefschetz type for singular 
non compact spaces, at the 
homotopy level. As applications, we derive results on the topology
 of the fibres of polynomial functions or 
the topology of complements of hypersurfaces in $\bC^n$.
   
\end{abstract}

\setcounter{section}{0}
\section{Introduction}

  In this paper we extend the method 
  of slicing by pencils to  
a larger class of ``admissible" pencils. This aim is motivated by the fact that,  
in certain situations, the pencils one is able to use are not generic. As an important example, we mention the study of the topology of a polynomial function $f \colon \bC^n \to \bC$, which is itself a nongeneric pencil on $\bC^n$ (see \S \ref{s:appl}).

 So let $X= Y\setminus V$, where $Y$ is a compact 
complex analytic space and $V$ is a closed complex subspace. For instance, quasi-projective varieties are of this kind.

 We call ``pencil" the ratio of two
sections $f$ and $g$ of a holomorphic line bundle $L\to Y$.
It defines a holomorphic function $h:= f/g$ over the complement $Y\m A$ of the ``axis" of the pencil (i.e. the indeterminacy locus) $A:= \{ f=g = 0\}$. A pencil is called {\em generic} with respect to $X$ when its axis $A$ is general (i.e. stratified
transversal to some Whitney stratification of the pair
$(Y,V)$)
and when the holomorphic map $h= f/g : Y\m A \to \bP^1$ has only stratified
double points as singularities. These singularities are finitely many,
by the compactness of $Y$, but part of those might be outside
$X$. Instead of only double points,
 one may consider pencils with any kind of stratified isolated singularities (see  \cite{HL}, \cite{GM}).
 
 We mean by {\em nongeneric pencil} a pencil where singularities may also occur inside the axis (e.g. when the axis $A$ is not transversal).
 In this paper we consider nongeneric pencils with
isolated singularities in the axis, in 
the precise meaning of Definition \ref{d:nongeneric}. 
The main result we prove is the following:

\begin{theorem}\label{th:main}
Let $X= Y\m V$, where $Y$ and $V\subset Y$ are compact complex analytic
spaces. 
Let $f, g$ be sections of a holomorphic line bundle over $Y$, defining a 
 pencil with at most isolated singularities in the axis (Definition 
\ref{d:nongeneric}) and let $X_\alpha$ denote a generic member of the pencil. Let $\rhd X \ge n$, where $n\ge 2$.

If one of the following two conditions is fulfilled:
\begin{enumerate}
\item  $A\not\subset V$ and the pair $(X_\alpha , A\cap X_\alpha)$ is $(n-2)$-connected,

\item $A\subset V$ and $V\supset \{ g=0\}$,

\end{enumerate}
then the pair $(X, X_\alpha)$  
is $n-1$ connected.
\end{theorem}
\begin{example}\label{e:main}
Let $\hat f = x^2y + xz^2$, $\hat g = z^3$ define a pencil on $\bP^3$ with homogeneous coordinates $x,y,z,w$. We restrict the pencil to the nonsingular surface $Y\subset \bP^3$ given by $yw+x^2 - z^2=0$. The axis $\hat A \subset \bP^3$ consists of two lines, one of which, namely $\{ x=z=0\}$ being the singular locus of each member of the pencil. The axis $A = \hat A\cap Y$ of the pencil on $Y$ is two points and one of them is singular. This is clearly a nongeneric pencil, since all of its members have singularities in the axis. 

Nevertheless, Theorem \ref{th:main} can be applied to this situation, since the pencil has isolated singularities in the axis (by Proposition \ref{l:lefsing}), $\rhd Y \ge 2$ and $(Y_\alpha , A\cap Y_\alpha)$ is $0$-connected.
\end{example}

Theorem \ref{th:main} represents a far reaching extension of the Lefschetz Hyperplane Theorem, which asserts that, if  $X\subset \bP^N$ is a projective variety and $H\subset \bP^N$ a hyperplane
 such that $X\setminus H$ is non singular of dimension $\ge 
n$ (more generally: if $X\setminus H$ has {\em rectified homotopy depth} greater or equal to $n$), then $\pi_k(X, X\cap H) =0$ for all $k<n$. Although not generic itself, the hyperplane $H$ can be viewed as a member of a generic pencil 
(by choosing a generic axis inside $H$, which is possible to do in the projective space). Then one may conclude by using our 
Remark \ref{c:nongen}, which deals with this particular situation.

Lefschetz's original proof (\cite{Lef}, see also \cite{La}) uses a generic pencil of hyperplanes to scan the space. Several generalizations to (non compact) spaces with singularities, such as by Goresky and MacPherson \cite{GM}, Hamm and L\^e [HL1-3], use Morse theory (method first employed by Bott \cite{Bo}, Andreotti and Frankel \cite{AF}). For the bibliography up to '88, one can look up \cite{GM}.
 It appears that in some of these generalizations, under the respective hypotheses, generic pencils do exist and their use yields alternative
 proofs (e.g. in \cite[Thm. 1.2, pag.199]{GM}).

The main condition we impose in our Theorem \ref{th:main} is on the rectified 
homotopical depth of the space $X$ (abbreviated $\rhd X$). This was introduced by Hamm and L\^e \cite{HL}, who proved several Lefschetz type theorems and Grothendieck's conjectures [HL-1,2,3,4]. We 
only say here that this amounts to a local condition $\rhd_x X \ge n$ which 
is satisfied, for instance, at points $x$ where $(X,x)$ is a germ of a complete intersection of dimension $n$. The condition on $\rhd$ is recursive, namely $\rhd X_\alpha \ge n-1$, by 
\cite[Theorem 3.2.1]{HL}.

The condition on the connectivity of $(X_\alpha, A\cap X_\alpha)$ comes in naturally (see \cite{La}) and can also be recursive. Indeed, $A\cap X_\alpha$ may be a (singular) fibre of a second pencil 
$h' = f/g'$ on the space $X_\alpha$ (which replaces $X$ as total space), 
with new axis $A_1 := \{ f=g' =0\}$. 
 Let's denote by 
$(X_\alpha)_\beta$ a generic member of it. We suppose by induction that the pair $((X_\alpha)_\beta, A_1\cap 
(X_\alpha)_\beta)$ is $(n-3)$ connected and suppose that Theorem \ref{th:main} can be applied. It follows that $(X_\alpha, 
(X_\alpha)_\beta)$ is $(n-2)$ connected. When $X$ is compact or when $A\cap X_\alpha$ is generic, this implies in turn, see Remark \ref{c:nongen}, that $(X_\alpha, A \cap X_\alpha)$ is $(n-2)$ connected.


 By applying Switzer's result \cite[Proposition 
6.13]{Sw} to the above theorem, one derives the usual
attaching result:

\begin{corollary}\label{c:main2}
 Under the hypotheses of Theorem \ref{th:main}, up to homotopy type, the space $X$ is built 
from $X_\alpha$ by attaching cells of dimension $\ge n$.
 
If $X$ 
is in addition a Stein space of dimension $n$, then the attaching cells are 
of dimension precisely $n$. If moreover $X$ is $n-1$ connected, then
the general hyperplane section $X_\alpha$  has the homotopy type of a 
bouquet of spheres $\vee S^{n-1}$. 
\fin
\end{corollary}


The proof of our theorem is based 
on the Nash blowing-up along the axis of the pencil, on homotopy excision and on 
local Lefschetz type results.
We also remark (Proposition \ref{p:gen1}) that, if $X_\alpha$ is not generic, 
then, in the conclusion of Theorem \ref{th:main} one may 
replace $X_\alpha$ by $X_D$, a small ``tube" neighbourhood of $X_\alpha$.

As natural applications, we prove new connectivity estimations on fibres of polynomial functions and complements of affine hypersurfaces.
\section{Singularities in the axis}

Fix a Whitney stratification $\cW$ of $Y$ such that $V$ is a
union of strata. 
Let $\tau = [s:t]$ denote a point on the complex projective line $\bP^1$. 
Let
$\bY$ denote the hypersurface $\{ (x,\tau) \in Y\times \bP^1 \mid
sf(x) - tg(x) =0\}$ in $Y\times \bP^1$ obtained by blowing up the axis
$A$ and let $\bX := \bY \cap (X
\times \bP^1)$. Consider the projection $p : \bY \to \bP^1$ to $\bP^1$
and its restriction $p_{|\bX} : \bX \to \bP^1$.
We also consider the projection to the first factor $\pi : \bY \to Y$.
 We use the following notations throughout the paper: for any
 $M\subset \bP^1$, $\bY_M := p^{-1}(M)$ and $\bX_M := \bX \cap
 \bY_M$.

Observe that
$A\times \bP^1 \subset \bY$, that
 $\bY\m(A\times \bP^1)$ can be identified (as the graph of $h$) with $Y\m A$, and that the restriction $p_{|\bY\m(A\times \bP^1)}$ can be identified 
 with $h$.
The stratification $\cW$ restricted to the open set $Y\m A$ induces a Whitney stratification on $\bY\m(A\times \bP^1)$, via the above identification.

\begin{definition}\label{d:str}({\bf Stratification of $\bY$})\\
Let $\cS$ be the coarsest Whitney
stratification on $\bY$ which coincides over $\bY\m (A\times\bP^1)$ with
the one induced by $\cW$ on $Y\m A$. (This exists, by classical
stratification arguments, see e.g. \cite{GLPW}.) We refer to it as the
canonical stratification of $\bY$ generated by the stratification $\cW$ of
$Y$. We also consider as canonical stratification of $\bX$ the stratification
induced by $\cS$ on $\bX$.
\end{definition}


After endowing $\bY$ and $\bX$ with canonical stratifications, the next
observation is that both $p: \bY \to \bP^1$ and $p_{|\bX}:
\bX \to \bP^1$ are {\em stratified
locally trivial fibrations above the complement in $\bP^1$ of some
finite
set} ("bad values") and that this is true in general, regardless of the
singularities or the position of the
axis $A$. This type of result is classical (Isotopy Theorem), it goes
back to Thom's paper \cite{Th} and it is based on the fact that $p$ is proper 
analytic and 
$\cS$ has finitely many strata. The problem one has to deal with is what
happens in the pencil when one encounters such a critical value. This comes
from the singular locus of $p$. Let us first define it.

\begin{definition}\label{d:sing}
 The {\em singular locus} of $p$ with respect to $\cS$ is the following closed
analytic subset of $\bY$:
 \[  \Sing_\cS p := \bigcup_{\cS_\beta \in \cS} \Sing
p_{|\cS_\beta}.\]
The {\em critical values} of $p$ with respect to $\cS$ are the points in the
image $p(\Sing_\cS p)$.
\end{definition}


\begin{definition}\label{d:nongeneric}
We say that the pencil defined by $h$ is a {\em (nongeneric) 
pencil 
with isolated singularities in the axis} if $\dim \Sing_\cS p \le 0$.
\end{definition}

 In case that $Y$ is projective, the condition $\dim \Sing_\cS p \le 0$ in the above definition is equivalent to the condition: the singularites of the 
function  $p$ at the blown-up axis $A\times \bP^1$ are at most isolated.
Indeed, in this case the singularities of $p$ outside the axis are also 
isolated, by the following reason. Suppose that $p$ has non-isolated singularities and consider one component $\cC$ of the singular locus $\Sing_\cS p$. Since $\cC$ is necessarily contained in a single fibre of $p$, call it $p^{-1}(b)$, it has to intersects the axis $A\times \{ b\}$, which is a hypersurface in $p^{-1}(b)$. Then, at all points of $\cC \cap (A\times \{ b\})$, we have $\dim \cC \ge 1$. This ends our argument.

 We shall assume for the remainder of this paper that the singular
locus of $p$ is of dimension $0$, hence it consists of a finite number
of points. This assumption is satisfied for instance in the following particular but very significant case.


\begin{proposition}\label{l:lefsing}
 Let $Y\subset \bP^N$ be a projective variety endowed with some Whitney stratification and let $\hat h=\hat f/\hat g$ define a pencil of hypersurfaces in $\bP^N$.  Let $B$ denote the set of points on $\hat A \cap Y$ where some member of the pencil is singular or where $\hat A$ is not transversal to $\cW$. If $\dim B\le 0$ and the singular points of
$h: Y\m A\to \bP^1$ with respect to $\cW$ are isolated then $p$ has isolated singularities.
\end{proposition}
\begin{proof}
On $\bY\m (A\times \bP^1)$,
  $p$ is just $h$ and its singularities are isolated. The notation $A$ stays for $\hat A \cap Y$, as usual. 
  
  Next, let us remark that $\bY = \bH \cap (Y\times \bP^1)$, where $\bH = \{ x\in \bP^N, [s:t]\in \bP^1 \mid s\hat f(x)- t\hat g(x) =0\}$. The singularities of $\bH$ are contained into $\Sigma_{\hat A} = (\hat A \times \bP^1) \cap \{ s\partial \hat f - t \partial \hat g = 0\}$, which is at most a collection of lines, by hypothesis. We endow $\bH$ with the coarsest Whitney stratification. It follows that $\bH$ is transversal to the strata $\cW \times \bP^1$ of $Y\times \bP^1$,
  except eventually along $B\times \bP^1$. By using that a transversal intersection of Whitney stratified sets is Whitney \cite[p. 19]{GLPW}), it follows that the canonical stratification $\cS$ of $\bY$ restricted to $(A\m B) \times \bP^1$ contains only strata which are products by $\bP^1$.  Hence the projection $p$ is transversal to these strata.  
  Finally, the stratification $\cS$ may distinguish at most a finite number of points, as point strata, out of each line from the collection of projective lines $B\times \bP^1$. Then $p$, being a projection, is still 
transversal to the complement of these points in $B\times \bP^1$.
\end{proof}

\section{Proof of the main result and some consequences}\label{s:proof}

Denote in the following $A' := A\cap X$. For any $M\subset \bP^1$,
let us denote $Y_M := \pi(p^{-1}(M))$
  and $X_M := X\cap Y_M$.
  
  We compute here the homotopy groups $\pi_j(X, X_\alpha)$ of the
pair space-section, where in the sequel $\alpha$ is supposed to be a 
general 
value for 
$p$.
Let $\Sing_\cS p =\{ b_1, \ldots , b_k\} \subset \bP^1$ be the
singular
 values of $p$.


\begin{lemma}\label{l:htdisc}
For any $W\subset \bP^1$, the space $Y_W$ (resp. $X_W$) is homotopy 
equivalent 
to the 
space $\bY_W$ (resp. $\bX_W$), to which one attaches along the product of 
$A$ (resp. $A'$) by $W$ 
the product of $A$ (resp. $A'$) by $\Cone(W)$.
In particular, if $W$ is contractible, then 
$Y_W \stackrel{\h}{\simeq} \bY_W$ and $X_W \stackrel{\h}{\simeq} \bX_W$.
\end{lemma}
\begin{proof}
The first statement immediately follows from the definition of the spaces
 $\bY$ and $\bX$.
As for the second statement, we have the homotopy equivalences:
\[ Y_W \stackrel{\h}{\simeq} \bY_W \cup_{A\times W} A\times \Cone(W) 
\stackrel{\h}{\simeq} \bY_W,\]
since $A\times \Cone(W)\stackrel{\h}{\simeq} A\times W$.
 The same argument applies to $X_W$ and we get $X_W \stackrel{\h}{\simeq} 
\bX_W$. Notice that in case $A' = \emptyset$, we have $\bX_W = X_W$, for 
any 
$W$.
\end{proof}

Take small disjoint closed discs $D_i\subset \bP^1$ centered at
$b_i$, a point $\a\in \bP^1$ exterior to all discs, and simple
paths (non self intersecting, mutually non intersecting except at $\alpha$) 
$\gamma_i
\subset \bP^1 \m \cup_{i=1}^{k} \stackrel{\circ}{D_i}$ from $\a$
to some fixed point $c_i\in \partial D_i$, for all $i\in
 \{ 1, \ldots , k\}$.
Denote by $\cD$ the subset $\bigcup_{i=1}^k D_i \cup \gamma_i$.
Choose a closed disc $K\subset \bP^1\m \stackrel{\circ}{\cD}$ such that 
$K\cap 
\cD = \{\alpha\}$ and consider the decomposition
 $\bP^1
= K \cup \overline{\bP^1 \m K}$, where $K\cap \overline{\bP^1 \m K}$ is a
 circle which we denote by $S$ (therefore $\alpha \in S$). 
 Since $\cD$ is a deformation retract of $\overline{\bP^1 \m K}$, 
and by the local triviality of the map $p_{|\bX}$, the space
 $\bX_{\overline{\bP^1 \m
K}}$ retracts to $\bX_\cD$. By Lemma \ref{l:htdisc}, the
corresponding retraction follows in the space $X$, namely
$X_{\overline{\bP^1 \m K}}$ retracts to $X_\cD$. Indeed:
\[ X_{\overline{\bP^1 \m K}}\stackrel{\h}{\simeq}\bX_{\overline{\bP^1 \m 
K}}
\cup_{A\times \overline{\bP^1 \m K}} A\times\Cone(\overline{\bP^1 \m K})
\stackrel{\h}{\simeq} \bX_\cD \cup_{A\times \cD} A\times\Cone(\cD) 
\stackrel{\h}{\simeq} X_\cD.\]
The same lemma
also shows that $X_\alpha$ is homotopy equivalent to $X_K$. Hence
we have the homotopy equivalence of pairs $(X, X_\alpha) \stackrel{\h}{\simeq} (X,X_K)$ and $(X_{\overline{\bP^1 \m K}},
X_\alpha) \stackrel{\h}{\simeq} (X_\cD, X_\alpha)$.

We now want to apply homotopy excision (Blakers-Massey theorem \cite{BM}, see also
\cite[Corollary 16.27]{Gr}) to the pair $(X_K \cup X_{\overline{\bP^1 \m K}}, X_K)$. We
have clearly that the pair $(X_{\overline{\bP^1 \m K}}, X_S)$ is $0$-connected and we
need the connectivity level of $(X_K, X_S)$. By considering the
triple $(X_K, X_S, X_\alpha)$ and remembering that $X_K
\stackrel{\h}{\simeq} X_\alpha$, we get, for any $i$, the
isomorphism:
\begin{equation}\label{eq:1}
\pi_{i+1}(X_K, X_S)\simeq \pi_{i}(X_S, X_\alpha).
\end{equation}

\begin{proposition}\label{p:bm}
 Assume that $A' \not =\emptyset$. 
\begin{enumerate} 
\rm \item \it If $(X_\alpha, A')$ is $m$-connected, $m\ge 
0$, then:
 \begin{enumerate}
\rm \item \it $(X_S, X_\alpha)$ is at least $m+1$ connected.
\rm\item \it  The excision morphism $\pi_j(X_{\overline{\bP^1 \m K}}, X_S) \to \pi_j(X,
X_K)$ is an isomorphism for $j\le m+1$ and an epimorphism for
$j=m+2$.
\end{enumerate}
\rm \item \it If, for any $i\in \{ 1,\ldots , k\}$, the pair $(\bX_{D_i} , 
\bX_{c_i})$ is $s$-connected, then $(X_\cD, X_\alpha)$ is $s$-connected 
too.
\end{enumerate}
\end{proposition}
\begin{proof}
 (a)(i).
 Note first that $X_S$ is homotopy equivalent to the subset
 $\bX_S \cup  A'\times K$ of $\bX_K$.
 Let $I$ and $J$ be two arcs (of angle less than $2\pi$) which cover
 $S$. We have the homotopy equivalence $(X_S, 
X_\alpha)\stackrel{\h}{\simeq}
 (\bX_I\cup  (A'\times K) \cup \bX_J \cup  (A'\times K)), \bX_J \cup
 ( A'\times K))$.

 Then, by homotopy excision (Blakers-Massey theorem), if the (identical!) 
pairs
 $(\bX_I\cup (A'\times K), \bX_{\partial I}\cup (A'\times K))$ and 
$(\bX_J\cup 
(A'\times K), \bX_{\partial 
J}\cup 
(A'\times K))$ are $m+1$
 connected, then the following morphism:
 \[ \pi_j(\bX_I\cup  (A'\times K), \bX_{\partial I}\cup  (A'\times K))
  \to \pi_j(\bX_I \cup  (A'\times K)\cup \bX_J \cup (A'\times K), \bX_J 
\cup
 (A'\times K))\]
is an isomorphism for $j\le 2m+1$. This would imply that $(X_S,
X_\alpha)$ is $m+1$ connected.

 It remains to prove our hypothesis. This follows since the
 pair $(\bX_I\cup  (A'\times K), \bX_{\partial I}\cup (A'\times K))$ is
 homotopy equivalent to $(\bX_\alpha\times I, \bX_\alpha\times\partial 
I\cup 
A'\times I)$
 and this, in turn, is just the product of pairs
 $(\bX_\alpha, A')\times (I, \partial I)$. Here we use the
 assumption in the statement of our Proposition.\\
 (ii). From (1) and from the point (i) it follows that the pair 
 $(X_K, X_S)$ is $m+2$ connected. Hence we've got the needed connectivity 
level 
of 
$(X_K, X_S)$ which makes homotopy excision work. The proof is now 
complete.\\
 (b).  From 
Lemma \ref{l:htdisc} it follows that $(X_\cD, X_\alpha) 
\stackrel{\h}{\simeq}(\bX_\cD, \bX_\alpha)$, since $\cD$ is contractible.
 By Switzer's result \cite[6.13]{Sw}, the $s$-connectivity of the 
CW-relative 
complex $(\bX_{D_i} , \bX_{c_i})$ implies that, up to homotopy equivalence, 
$\bX_{D_i}$ is obtained from  $\bX_{c_i}$ by attaching cells of dimension 
$\ge 
s+1$. Since $\bX_\cD = \cup_i \bX_{D_i \cup \gamma_i}$ and 
 $\bX_{c_i}\stackrel{\h}{\simeq} \bX_\alpha$, it follows that $\bX_\cD$ is 
obtained from  $\bX_\alpha$ by attaching cells of dimension $\ge s+1$.
\end{proof}


\begin{proposition}\label{c:decomp}{\bf (case $V\supset \{ g=0\}$)}\\
Assume that $V\supset \{ g=0\}$. If the pair $(\bX_{D_i} , \bX_{c_i})$ is 
$s$-connected, for any $i\in \{ 1,\ldots , k\}$,  then $(X, X_\alpha)$ is 
$s$-connected.
\end{proposition}
\begin{proof}
 In this case $h_{|X} = f/g :X\to \bC$ is a well defined holomorphic 
function. 
Since $\bC$ retracts to $\cD$, the spaces $K$ and $S$ do not occur and we 
simply have $(X, X_\alpha)\stackrel{\h}{\simeq}(X_\cD, X_\alpha)$. 
 We have $A' = \emptyset$, therefore $\bX_{D_i} = X_{D_i}$ and $\bX_{c_i} 
\stackrel{\diff}{\simeq} X_\alpha$. The result follows by using Switzer's 
argument \cite[6.13]{Sw}, like in the proof of Proposition \ref{p:bm}(b).
\end{proof}

Since $A\subset V$, the axis $A$ may be highly non transversal to the 
stratification $\cW$ of $(Y,V)$. However, the singularities of $p$ might be still 
isolated; remember that Definition \ref{d:sing} says that 
$x\in \Sing_\cS p$ if and only if $x\in (\Sing Y_t) \cap A$, where $t= 
p(x)$.
This situation occurs when studying polynomial functions on $\bC^n$.
For instance, the polynomial function $f\colon \bC^2 \to \bC$, $f(x,y) = x+x^2y$, as pencil of hypersurfaces, has isolated singularities at infinity. See \S \ref{s:appl} for applications.

\vspace*{5mm}

 We have shown that in the cases covered by Propositions \ref{p:bm} and \ref{c:decomp}, the connectivity of 
$(X, X_\alpha)$ depends on the one of $(\bX_{D_j}, \bX_{c_j})$, for each $j$. So we further study $(\bX_{D_j}, \bX_{c_j})$, for some fixed $j$. Let's drop the index $j$ and write
 simply $(\bX_{D}, \bX_{c})$. We have assumed that $p$ has no singularities 
over $D^*$ and has only isolated singularities over the center $b$ of $D$, 
among the following three possible kinds: singularities on
  $\bX_{b} \m A\times \bP^1
  $, singularities on
  $\bY_{b} \m (\bX \cup A\times \bP^1)$ and
  singularities of $p$ in the axis $A\times \{ b\}$.

Say $\bY_b \cap \Sing_\cS p = \{ a_1, \ldots , a_r\}$. Let us consider 
(small) 
local Milnor-L\^e balls at each isolated singularity of $\bY_b$. The 
existence 
of such balls was shown by Milnor \cite{Mi} in case of smooth ambient space 
and 
by L\^e D.T. for singular stratified spaces \cite{Le-oslo}, \cite{Le-jag}. 
These 
closed balls have the property that, modulo the reducing of the radius of 
$D$ 
as much as necessary, their boundaries $\partial B_i$
are transversal to the strata of our stratification $\cS$ of the space 
$\bY$.
  It is natural now to excise the complement $C$ of the disjoint union 
$\sqcup_{i\in \{ 1, \ldots , r\}} B_i \cap \bY_\cD$ from the pair $(\bY_D, 
\bY_c)$. This has to be related to the fact that $p_| : C\to D$ is a 
trivial 
fibration, and moreover, since it is a stratified fibration, the 
restriction 
$p_{|\bX} : C\cap \bX \to D$ is also trivial.
  
  If we perform excision, then we reduce the problem to a local one, around 
the 
isolated singularities.  In homology, we get the direct sum  decomposition 
$H_* 
(\bX_D, \bX_c) = \oplus_i H_* (B_i \cap\bX_D, B_i \cap\bX_c)$. In homotopy, 
the 
excision (Blakers-Massey theorem) implies that the level of connectivity of 
$(\bX_D, \bX_c)$ is at least equal to the minimum of the levels of 
connectivity 
of  $(B_i \cap\bX_D, B_i \cap\bX_c)$.
  
  We need a condition which implies a certain level of connectivity of 
each 
pair $(B_i \cap\bX_D, B_i \cap\bX_c)$. A condition that fits well is the  
rectified homotopical depth of the total space $\bX$. This condition does 
not 
depend on the stratification of the space.
  

\begin{proposition}\label{p:rhd}
If $\rhd(X) \ge s+1$, then, for any $i\in \{ 1, \ldots , r\}$, the pair 
$(B_i 
\cap\bX_D, B_i \cap\bX_c)$ is at least $s$-connected.
\end{proposition}
\begin{proof}
 Since we work on the space $\bX$, we would need a condition on 
$\rhd(\bX)$.
 So, we first prove that $\rhd(\bX)\ge s+1$. Since on the space $X\times 
\bP^1$ 
we have the product stratification $\cW\times \bP^1$, the condition 
$\rhd_\cW(X) \ge s+1$ implies $\rhd_{\cW\times \bP^1}(X\times\bP^1)\ge 
s+2$. 
Our space $\bX$ is a hypersurface in  $X\times\bP^1$ and therefore its 
rectified homotopical depth is one less, i.e. $\rhd\bX \ge s+1$, by 
\cite[Theorem 3.2.1]{HL}.
 
 The rectified homotopical depth of $\bX$ gives a certain level of 
connectivity 
of the complex links of the strata of the stratification $\cS$, according 
to 
\cite{GM} and \cite{HL}. One may relate the connectivity of these complex 
links to the 
connectivity of the Milnor-L\^e data $(B_i \cap\bX_D, B_i \cap\bX_c)$. This 
is 
more special data, especially when the singularity is not in $\bX_D$ but on its 
``boundary" $\bY \cap (V\times \bP^1)$. Such relation among connectivities 
is 
the local Lefschetz theorem of Hamm and L\^e 
\cite[Theorem 4.2.1 and Cor. 4.2.2]{HL}. This result can be applied for the 
function $p_| : \bY_D \to D$ with isolated singularities and for the space 
$\bX_D 
= \bY_D \m (V\times \bP^1)$ and it tells precisely that, since $\rhd \bX 
\ge 
s+1$, the pair $(B_i \cap \bX_D, B_i\cap \bX_c)$ is at least $s$-connected.
This proves our statement.
\end{proof}   


\begin{proof} {\it of Theorem \ref{th:main}.} 
From the long exact sequence of the triple $(X_{\overline{\bP^1\setminus K}}, X_S, X_\alpha)$ and 
since $(X_S, X_\alpha)$ is $n-1$ connected (by Proposition 
\ref{p:bm}(a)(i)), 
it follows that the morphism 
(induced by inclusion):
\begin{equation}\label{eq:2}
\pi_j(X_{\overline{\bP^1\setminus K}},X_\alpha) \to \pi_j(X_{\overline{\bP^1\setminus K}},X_S) 
\end{equation}
 is an isomorphism for $j\le n-1$ and an epimorphism for $j= n$.
 
  Next, by Proposition \ref{p:rhd} and the observations before it, the 
pairs 
$(\bX_{D_i}, \bX_{c_i})$ are $n-1$ connected. Then, applying Proposition 
\ref{p:bm}(b) and (a)(i), via the morphism (\ref{eq:2}), we get the 
connectivity 
 $n-1$ of $(X, X_\alpha)$. So far for the proof of Theorem \ref{th:main}(a).
 
 As for (b), it is now an imediate consequence of 
Proposition \ref{c:decomp} together with Proposition \ref{p:rhd}.
Note however something important, which we shall use in \S \ref{s:appl}: the hypothesis $\dim \Sing_\cS p \le 0$ is too strong. Actually, in Theorem \ref{th:main}(b) it 
is sufficient to assume that $\dim \Sing_{\cS} p' \le 0$, where $p'$ is the restriction of $p$ on $\bY_{\bC}$. Our proof only uses this weaker 
hypothesis.
\end{proof}   

We end this section by proving the usual observation that we can replace $X_\alpha$ in the 
conclusion of the theorem by a ``very bad'' member of the pencil 
and then take a good neighbourhood of this. 
Let $D_\delta \subset \bP^1$ denote a small enough disc centered at 
$\delta$ such that $D_\delta  \cap p(\Sing_\cS p) = \{ \delta\}$.

\begin{proposition}\label{p:gen1}
In Theorem \ref{th:main}, replace the 
hypothesis about the singularities of the pencil with the following: ``Let 
$f,g$ define a pencil with at most isolated singularities except at one 
fiber $X_\delta$, that is $\dim (\Sing_\cS p \cap X_\beta )\le 0$, for all 
$\beta \not= \delta$''.

Then $(X, X_{D_\delta})$ is $(n-1)$ connected.
\end{proposition}
\begin{proof}
 We just consider 
$D_\delta$ as one of the small discs $D_i$ within $\cD$. We of course still need a generic member $X_\alpha$ and the hypothesis on the connectivity of 
$(X_\alpha, A')$ in case $A'$ is not empty. 

It follows from Proposition \ref{p:bm}(a)(i) and 
from (\ref{eq:1}) that $(X_K, 
X_S)$ is $n$-connected. By homotopy excision of $X_\cD$ from $(X, X_\cD)$, 
this 
implies that $(X, X_\cD)$ is $n$-connected. In turn, via the homotopy exact 
sequence of the triple $(X, X_\cD, X_{D_\delta})$, this implies that the 
morphism 
induced by inclusion $\pi_i (X_\cD, X_{D_\delta}) \to \pi_i (X, 
X_{D_\delta})$ is an 
isomorphism for $i\le n-1$. It remains to study $(X_\cD, X_{D_\delta})$. 
Since $\cD$ and ${D_\delta}$ are 
contractible, Lemma \ref{l:htdisc} says that $(X_\cD, X_{D_\delta}) 
\stackrel{\h}{\simeq} (\bX_\cD, \bX_{D_\delta})$. 
By the argument in the proof of Proposition \ref{p:bm}(b) and also by 
Proposition \ref{p:rhd}, $\bX_\cD$ is obtained from $\bX_{D_\delta}$ by 
attaching cells of dimension $\ge n$.
Consequently, $(\bX_\cD, \bX_{D_\delta})$ is also 
$(n-1)$-connected, hence $(X, X_{D_\delta})$ too.
\end{proof}   

Suppose now that $X$ is compact. Then, since in generic pencils on such $X$ the tube neighbourhood $X_{D_\delta}$ contracts to
the fibre $X_{\delta}$, we have:  

\begin{remark}\label{c:nongen}
 In the conclusion of the above Proposition \ref{p:gen1}, the tube $X_{D_\delta}$ can be replaced by $X_\delta$ when the inclusion $X_\delta \subset X_{D_\delta}$ is a homotopy equivalence. This happens for instance if $X$ is compact or if the  pencil has no singularities in the 
axis (i.e. $(A\times\bP^1) \cap \Sing_\cS p = \emptyset$).
\end{remark}

 As mentioned in the Introduction, this remark is useful in order to make 
Theorem \ref{th:main} 
work inductively on dimension, namely to show that the condition on the 
connectivity of $(X_\alpha , A\cap X_\alpha)$ is recursive. It also clarifies 
why the ``classical'' Lefschetz hyperplane theorem (stated in the 
Introduction) can be proved by using a pencil having a generic axis into the ``bad'' hyperplane $H$.

\section{Aplications}\label{s:appl}
We give applications of Theorem 
\ref{th:main}(b) concerning the topology of fibres of polynomial functions $f: \bC^n 
\to \bC$ and complements of hypersurfaces in $\bC^n$.
The asymptotic behaviour of $f$ was 
studied 
by several authors in the last years
(see e.g. \cite{ST}, \cite{Pa}, \cite{Ti-reg}). 

  A polynomial function may be extended to a meromorphic function on a compact space; this embedding is however not unique. We consider 
here the embedding of $\bC^n$ into a weighted projective space $\bP_w := \bP(w_1, \ldots , w_n, 1)$, as follows. 
Associate to each coordinate $x_i$ a positive 
weight $w_i$, for $i\in \{ 1, \ldots , n\}$ and write $f = f_d + f_{d-k} + 
\cdots$ where $f_i$ is the degree $i$ weighted-homogeneous part of $f$ and 
$f_{d-k} \not= 0$.  Then take a new variable $z$ of weight 1. We get a 
meromorphic function $\tilde f/z^d$ on $\bP_w$ and the Nash blown up space 
$\bY := 
\{ s\tilde f - tz^d =0\} \subset \bP_w \times \bP^1$, where  
$\tilde f$ is the degree $d$ homogenized of $f$.  
We consider the coarsest Whitney stratification $\cS$ on $\bY$.
  Here $V$ is the hyperplane at 
infinity $\{ z=0\}$ and $g := z^d$.

We say that the pencil defined by $f$ 
{\em has at most isolated singularities in $\bY' := \bY_{\bC}$} if the 
singularities of the restriction $p':=p_{|\bY'} : \bY' \to \bC$ with respect 
to the 
stratification $\cS$ are isolated.
Let us denote by $\Sigma$ 
the weighted projective variety $\{ \grad f_d =0, f_{d-k}=0\} \subset 
\bP_w\cap 
\{ z=0\}$.


\begin{proposition}\label{l:weight}
If $\dim \Sing f\le 0$ and $\dim \Sigma \le 0$ then the pencil defined by $f$ 
has at most isolated singularities in $\bY'$. 
\end{proposition}
\begin{proof}
 We prove that $\dim \Sing_\cS p'\le 0$. Since 
the singularities of $p'$ on $\bC^n = \bY' \setminus \{ z=0\}$ are isolated 
by 
hypothesis, we only have to look on $\bY' \cap \{ z=0\}$, which is the 
product
$\{f_d =0\} \times \bC \subset (\bP_w \cap \{ z=0\})\times \bC$.

We need to know the stratified structure of $\bY'$ in the neighbourhood of 
$\bY' \cap \{ z=0\}$.
 Consider first $\bY' \cap \{ z=0\}\cap (\{ \grad f_d \not= 0\} \times \bC)$. 
This 
is a subspace of the quotient by the $\bC^*$ action of the nonsingular part 
of
the hypersurface $\{ \tilde f - tz^d =0\}$ considered as subset of 
$(\bC^{n+1} \setminus \{ 0\})\times \bC$. We claim that the stratification by 
the orbit 
type of this nonsingular part is Whitney regular. Indeed, the $\bC^*$ action 
reduces, within Zariski-open subsets $\{ x_i\not= 0\}$, to the action of a 
finite group. For quotients by finite 
groups, the natural orbit type stratification is Whitney regular (see e.g.
 \cite[p. 21]{GLPW},
\cite{Fe}). 
 
Since the action on the factor $\bC$ is trivial, the 
strata within $\bY' \cap \{ z=0\}\cap (\{ \grad f_d \not= 0\} \times \bC)$ 
are products by 
$\bC$. It follows that $p'$ (which is the projection to $\bC$) is 
transversal to these strata.

Next we look to points 
$q\in \{ \grad f_d =0\} \setminus \{ f_{d-k} = 0\} \subset \bP_w \cap \{ 
z=0\}$.   We have $\tilde f(x,z) - tz^d = f_d(x) + z^k g$, 
where $g(q,z,t) = f_{d-k}(q) + z h(q,z,t)$ 
and $f_{d-k}(q) \not= 0$.
 
 Locally at $q$, our hypersurface $\bY' = \{ \tilde f (x,z) - tz^d 
=0\}$ is equivalent, via an analytic change of 
coordinates, modulo a choice of the $k$-root, to the product of 
 $\{f_d(x) + (z')^k = 0\}$ by the coordinate $t$, where $z'$ is the new 
coordinate $z$ and $\{z' =0\}$ is equal to $\{z =0\}$ at $q$.
  Since we have again a product by $\bC$, we may deduce that, locally, the 
line $\{ q\} \times \bC$ is included into a Whitney 
stratum of $\cS$. Hence the map $p'$ has no singularities on $\{ q\} \times 
\bC$.

   This proves that $\Sing_\cS p' \subset \Sigma \times \bC$. Since 
$\dim \Sigma \le 0$, the set $\Sigma \times \bC$ is a finite union of complex 
lines. The map $p'$ is transversal to such a line, so singularities of $p'$ 
on 
$\Sigma \times \bC$ can occur only if $\Sigma \times \bC$ contains 
point-strata 
of $\cS$. But there can only be finitely many such point-strata. This ends 
our 
proof.
\end{proof}

\begin{note}\label{n:dp}
The above proof shows that the condition $\dim \Sigma \le 0$ implies that $\dim \Sing f \le 1$.
\end{note}

\begin{corollary}\label{p:weight}  
If the polynomial $f: \bC^n \to \bC$ has isolated singularities and $\dim 
\Sigma \le 0$ then its general fibre $X_\alpha$ is homotopy equivalent to a bouquet of spheres of dimension $n-1$. Moreover, any atypical fibre of $f$ is at least $n-3$ connected.
\end{corollary}
\begin{proof}
 By applying Theorem \ref{th:main}(b) via Proposition \ref{l:weight}, we get that $(X, 
X_\alpha) = (\bC^n, X_\alpha)$ is $n-1$ connected. Since $X_\alpha$ is Stein of 
dimension $n-1$, our first statement follows by Corollary \ref{c:main2}.

For the second statement, we use an argument similar to the one used to prove \cite[Theorem 5.5]{Ti-compo}. We give here an outline of it.
First, we take a sufficiently large ball $B$ centered at the origin such that
$X_b \cap B_M$ is homotopy equivalent to $X_b$, for any ball $B_M$ larger or equal to $B$.

We next consider the polar locus of the map $(p',z)$ at some isolated
singularity in $\bY_b \cap \{ z=0\}$. This locus is a curve or it is void, see
{\em loc. cit.} The polar curves corresponding to the singularities 
in $\bY_b \cap \{ z=0\}$ intersect a nearby general fibre $X_\alpha$ at a finite number of points. It follows that $X_\alpha$ is homotopy equivalent to $X_\alpha \cap B$ to which one attaches a finite number of  cells of dimension $n-1$, corresponding to the intersection multiplicities of the polar curves with $X_\alpha$. On the other hand, $X_b \cap B$ is homotopy equivalent to $X_\alpha \cap B$ to which one attaches $n$-cells, corresponding to the isolated singularities on $X_b$. It follows that $(X_\alpha, X_\alpha \cap B)$
is $n-2$ connected and that $(X_b, X_\alpha \cap B)$
is $n-1$ connected. Therefore, since $X_\alpha$ is $n-2$ connected, $X_b$ is $n-3$ connected.
\end{proof}

This represents an extension of results on connectivity of fibres,
in the vein of \cite{ST} and \cite{Ti-compo}. 
Dimca and P\u aunescu \cite{DP} proved recently a related result, by different methods. We shall explain in a forthcoming paper \cite{LT} that, with a recursive procedure, our result not only recovers \cite{DP}, but also improves the connectivity estimation for generic fibres.


\begin{example}\label{e:wt}
$f = x^2y^3 + v^2 + y$ has $\Sing f = \emptyset$ and $\Sigma = \{ [1:0:0], 
[0:1:0]\} \subset \bP^2$. Here all the weights are $1$ and $\bP_w$ is the 
usual 
projective space $\bP^3$. One can verify that the general fibre is homotopy 
equivalent to $S^2\vee S^2\vee S^2$. It turns out that the only atypical fibre is $f^{-1}(0)$.
\end{example}


 We further give an application of our results to the topology of complements
 of hypersurfaces in $\bC^n$.
  This is a topic which goes back to Zariski and van Kampen \cite{vK}, who 
described a general procedure to compute the fundamental group of the 
complement to an algebraic curve in $\bP^2$ by slicing with linear pencils.
   Zariski showed that $\pi_1$ depends on the type and position of 
singularities of the curve. More recently, Libgober \cite{Li} proved similar 
results on the higher homotopy group $\pi_{n-k-1}(\bC^n \setminus V)$, where 
$k$ is the dimension of the singular locus of the hypersurface $V$, $k< n-2$. 
This is the first possible nontrivial homotopy group of rank higher than $1$, 
since, by the classical Lefschetz theorem, $\pi_{i}(\bC^n \setminus V) = 0$, 
for $2\le i < n-k-1$. Since $\pi_{n-k-1}(\bC^n \setminus V) =  
\pi_{n-k-1}(H_k \cap \bC^n \setminus V)$, for a general linear subspace $H_k$ 
of codimension $k$, the problem of finding $\pi_{n-k-1}$ reduces to the case 
of the complement of a hypersurface $V$ with isolated singularities (see 
\cite[\S 1]{Li}). 
   
   In \cite[Theorem 2.4]{Li}, Libgober considers the case when $V$ has at 
most isolated singularities and $\bar V$ is transversal to the hyperplane at 
infinity in $\bP^n$. We show here that, under certain circumstances (weaker 
transversality condition but imposing that $V$ is a generic fibre of $f$), we 
may conclude to the triviality of $\pi_{n-k-1}(\bC^n \setminus V)$.
 For a more general statement and other results on complements of hypersurfaces we refer to the forthcoming paper \cite{LT}. 

\begin{proposition} \label{p:vkampen}
If $V$ is a general fibre of a polynomial $f \colon \bC^n \to \bC$ and the pencil defined by $f$ has at most isolated singularities in $\bY'$ then $\bC^n \setminus V \stackrel{\h}{\simeq} 
S^1 \vee\bigvee S^n$. In particular, $\pi_{n-1}(\bC^n \setminus V) = 0$ if $n>2$.
\end{proposition}
\begin{proof}
We consider small enough discs $D_i \subset \bC$ centered at the ``bad 
values" of the pencil, like in \S 3, and also a small enough closed disc $D_0$ 
centered at the general value $\beta$. Let $\alpha\in \partial  D_0$. 
Using the notations in \S 3, the configuration $D_0^* \cup_i (D_i\cup 
\gamma_i)$ is a strong deformation retract of $\bC\setminus \{ \beta\}$. It 
follows that $\bC^n \setminus V \stackrel{\h}{\simeq} f^{-1}(D_0^* \cup_i 
(D_i\cup \gamma_i))$.

Now, denoting by $S^1$ the boundary of $D_0$, the space $f^{-1}(D_0^* 
\cup_i (D_i\cup \gamma_i))$ is homotopy equivalent to the space obtained by attaching to $f^{-1}(D_0^*) \stackrel{\h}{\simeq} S^1\times 
f^{-1}(\alpha)$ the space $f^{-1}(\cup_i(D_i\cup \gamma_i))$,
over $ \{\alpha\} \times f^{-1}(\alpha)$. Notice that  
$f^{-1}(\cup_i(D_i\cup \gamma_i))$
 is homotopy equivalent to $\bC^n$, thus attaching it to $S^1\times 
f^{-1}(\alpha)$ over $\{\alpha\} \times f^{-1}(\alpha)$ corresponds to shrinking $\{\alpha\} \times f^{-1}(\alpha)$ to a point.

We therefore get that $f^{-1}(D_0^* \cup_i 
(D_i\cup \gamma_i))$ is homotopy equivalent to a bouquet of
a circle $S^1$ and the suspension over $f^{-1}(\alpha)$.
Since by Corollary \ref{c:main2} and the remark following it, 
$f^{-1}(\alpha)$ is homotopy equivalent to a bouquet of $(n-1)$-spheres, the suspension over $f^{-1}(\alpha)$ is a bouquet of $n$-spheres. 
This concludes the proof.
 
\end{proof}

  In the Example \ref{e:wt}, the ``singularities at infinity" of a general fibre 
$V$, in the sense of \cite{Li} are non isolated (since $\bar V \not\pitchfork
  \{ z=0\}$ along $\{ x=0\} \cup \{ y=0\} \subset \bP^2$). The results of \cite{Li} do not give any information in this case; hoewever, our Proposition \ref{p:vkampen} can be applied, since the pencil 
defined by $f$ has at most isolated singularities in $\bY'$. We get $\bC^3 
\setminus V \stackrel{\h}{\simeq} S^1\vee S^3\vee S^3\vee S^3$.



\begin{thebibliography}{MMM}
\footnotesize{


\bibitem[AF]{AF}
A. Andreotti, T. Frankel, {\em
The second Lefschetz theorem on hyperplane sections}, in: Global Analysis,
Papers in Honor of K. Kodaira, Princeton Univ. Press 1969, p.1-20. 


\bibitem[BM]{BM}
A.L. Blakers, W.S. Massey, {\em
The homotopy groups of a triad. III},
Ann. of Math. (2)  {\bf 58} (1953), 409-417.

\bibitem[Bo]{Bo}
R. Bott,  {\em On a theorem of Lefschetz}, Michigan Math. J. {\bf 6} 
(1959), 211-216. 


\bibitem[DP]{DP}
A. Dimca, L. P\u aunescu, {\em On the connectivity of complex
affine hypersurfaces. II}, Topology {\bf 39}, 5 (2000),  1035-1043.

\bibitem[Fe]{Fe} 
M. Ferrarotti, {\em $G$-manifolds and stratifications}, Rend. Istit. Mat. 
Univ. Trieste {\bf 26}, no. 1-2 (1994), 211-232.

\bibitem[GLPW]{GLPW}
C.G. Gibson, K. Wirthm\"uller, A.A. du Plessis, E.J.N. Looijenga, {\em
Topological
Stability of Smooth Mappings}, Lect. Notes in Math. {\bf 552}, Springer
Verlag 1976.


\bibitem[GM]{GM}
M. Goresky, R. MacPherson,
{\em Stratified Morse theory},
Ergebnisse der Mathematik und ihrer Grenzgebiete. 3. Folge, Bd. 14.
Berlin Springer-Verlag 1988.

\bibitem[Gr]{Gr}
B. Gray,
{\em Homotopy theory. An introduction to algebraic topology},
Pure and Applied Mathematics. Vol. 64,  Academic Press 1975.

\bibitem[HL1]{HL-1}
H.A. Hamm,  L\^e D.T., {\em Lefschetz theorems on quasi-projective 
varieties},
Bull. Soc. Math. France {\bf 113} (1985), 123-142.

\bibitem[HL2]{HL-2}
H.A. Hamm,  L\^e D.T., {\em Local generalizations of Lefschetz-Zariski 
theorems}, J. Reine Angew. Math. {\bf 389} (1988), 157-189.

\bibitem[HL3]{HL}
H.A. Hamm,  L\^e D.T.,
{\em Rectified homotopical depth and Grothendieck conjectures},
The Grothendieck Festschrift, Collect. Artic. in Honor of the 60th
Birthday
of A. Grothendieck. Vol. II, Prog. Math. 87 (1990), 311-351.

\bibitem[HL4]{HL-4}
H.A. Hamm,  L\^e D.T., {\em Relative homotopical depth and a conjecture of 
Deligne}, Math. Ann. {\bf 296}, 1 (1993), 87-101. 

\bibitem[La]{La}
 K. Lamotke,
{\em The topology of complex projective varieties after S. Lefschetz}, 
Topology {\bf 20} (1981), 15-51. 

\bibitem[Lef]{Lef} 
S. Lefschetz,
{\em L'analysis situs et la geometrie algebrique}, Gauthier-Villars, Paris 
1924, nouveau tirage 1950.

\bibitem[L\^ e1]{Le-oslo}
L\^ e D.T., {\em Some remarks on the relative monodromy}, Real and Complex 
Singularities Oslo 1976, Sijhoff en Nordhoff, Alphen a.d. Rijn 1977, pp. 
397-403.

\bibitem[L\^ e2]{Le-jag}
 L\^ e D.T., {\em Complex Analytic Functions with Isolated
Singularities }, J. Algebraic Geometry {\bf 1} (1992),  83-100.

\bibitem[Li]{Li}
A. Libgober, 
{\em Homotopy groups of the complements to singular hypersurfaces. II}, 
Ann. Math., II. Ser. {\bf 139}, No.1 (1994), 117-144.

\bibitem[LT]{LT}
A. Libgober, M. Tib\u ar,
{\em Homotopy groups of complements and non isolated singularities}, manuscript.

\bibitem[Mi]{Mi}
J. Milnor, {\em Singular points of complex hypersurfaces}, Ann. of Math. 
Studies~61, Princeton 1968.

\bibitem[Pa]{Pa}
A. Parusi\'nski,  {\em A note on singularities at infinity of complex polynomials},
in: ``Simplectic singularities and geometry of gauge fields", Banach Center Publ.vol. {\bf 39} (1997), 131-141.


   \bibitem[ST]{ST}
 D. Siersma, M. Tib\u ar,  {\em Singularities at infinity and their
vanishing cycles},
 Duke Math. Journal {\bf 80}:3 (1995), 771-783.
 
\bibitem[Sw]{Sw}
R. Switzer, {\em Algebraic Topology - Homotopy and Homology}, Springer 
Verlag,
Berlin-Heidelberg-New York.

 \bibitem[Th]{Th}
R. Thom, {\em Ensembles et morphismes stratifi\' es}, Bull. Amer.
Math.
Soc. {\bf 75} (1969), 249-312.


 \bibitem[Ti1]{Ti-compo}
M. Tib\u ar, {\em Topology at infinity of polynomial maps and Thom
regularity condition}, Compositio Math. {\bf 111}, 1 (1998), 89-109.

\bibitem[Ti2]{Ti-reg}
M. Tib\u ar, {\em Regularity at infinity of real and complex polynomial 
maps},   
 in: Singularity Theory, The C.T.C Wall Anniversary Volume,  
LMS Lecture Notes Series {\bf 263} (1999), 249-264. Cambridge 
University Press.

\bibitem[vK]{vK}
E.R. van Kampen, 
{\em On the fundamental group of an algebraic curve},
Amer. J. Math. {\bf 55} (1933), 255-260.

 
 }
\end{thebibliography}
\end{document}